# PEIRCE'S TRUTH-FUNCTIONAL ANALYSIS AND THE ORIGIN OF TRUTH TABLES

Irving H. Anellis

*Abstract*. We explore the technical details and historical evolution of Charles Peirce's articulation of a truth table in 1893, against the background of his investigation into the truth-functional analysis of propositions involving implication. In 1997, John Shosky discovered, on the verso of a page of the typed transcript of Bertrand Russell's 1912 lecture on "The Philosophy of Logical Atomism" truth table matrices. The matrix for negation is Russell's, alongside of which is the matrix for material implication in the hand of Ludwig Wittgenstein. It is shown that an unpublished manuscript identified as composed by Peirce in 1893 includes a truth table matrix that is equivalent to the matrix for material implication discovered by John Shosky. An unpublished manuscript by Peirce identified as having been composed in 1883-84 in connection with the composition of Peirce's "On the Algebra of Logic: A Contribution to the Philosophy of Notation" that appeared in the *American Journal of Mathematics* in 1885 includes an example of an indirect truth table for the conditional.

**1. Introduction.** Charles Sanders Peirce (1839–1914) undertook a study of the conditions for the truth of propositions that continued through his entire career as a logician. The centerpiece of his work was the analysis of illation, equivalent when dealing with the propositional calculus with classical material implication. Although his notation and approach to logic shifted from syllogistic to algebraic to graphical, his interest in determining the conditions under which a concludion is true , and his preference throughout his work for the graphical, rooted in his study of matrix theory for linear and multilinear algebras, combined with his interest in establishing the conditions under which logical arguments lead to true conclusions, led him to devise a truth table matrix in 1893, thus nearly two decades ahead of the truth tables discovered by John Shosky [*1997*] attributable to Bertrand Russell and Ludwig Wittgenstein jointly, and even prior to those developed by Peirce for triadic logic in the period 1902-09 that were first announced by Atwell Rufus Turquette and Max Harold Fisch (1901–1995) ([Turquette *1964*] and [Fisch & Turquette *1966*]) and rediscovered by William Glenn Clark (1915–1993) and Shea Zellweger ([Clark *1997*] and [Zellweger *1997*]).

Shosky emphasized the difference between what he called the "truth table technique" and the "truth table device", meaning by the former truth-functional analysis of propositions and arguments, and the truth table or matrix array for displaying the truth values of propositions and arguments. In reply to Shosky [*1997*], [Anellis *2004*] explored the historical background to Peirce's development of truth table matrices in the work of William Stanley Jevons (1835–1882) [*1874*, 135] and Peirce's student Christine Ladd-Franklin (1847–1930) [*1883*] and argued that the history of the development of the truth table matrix is more complex than suggested by Shosky and that the earliest known example of a truth table matrix, in its familiar configuration whose authorship can be identified is not, contrary to Shosky's claim, that of 1912 ascribed to Russell, but that, for trivalent logic, presented by Peirce in a manuscript of 1902.

In the present account, I note the occurrence of a truth table matrix in a manuscript of Peirce of 1893 and examine the historical and conceptual context in which it occurs, against the background of Peirce's work on the truth-functional analysis of propositions. It transpires that the matrix presented by Peirce in 1893 is the table for an implicational inference between two terms, and is exactly equivalent to Russell's 1912 table for material implication.

**2. Peirce's development of the truth table matrix in the context of his study of implication.** In an undated, untitled, two-page manuscript designated "Dyadic Value System" (listed in the Robin catalog as MS #6; [Peirce *n.d.(a)*]), Peirce asserts that the simplest of value systems serves as the foundation for mathematics and, indeed, for all reasoning, because the purpose of reasoning is to establish the truth or falsity of our beliefs, and the relationship between truth and falsity is precisely that of a dyadic value system, writing specifically that "the the whole relationship between the values," 0 and 1 in what he calls a cyclical system "may be summed up in two propositions, first, that there are different values" and "second, that there is no third value." He goes on to says that: "With this simplest of all value-systems mathematics must begin. Nay, all reasoning must and does begin with it. For to reason is to consider whether ideas are true or false." At the end of the first page and the beginning of the second, he mentions the principles of Contradiction and Excluded Middle as central. In a fragmented manuscript on "Reason's Rules" of *circa* 1902 [Peirce *ca. 1902*], he examines how truth and falsehood relate to propositions.

Consider the formula $[(\sim c \supset a) \supset (\sim a \supset c)] \supset \{(\sim c \supset a) \supset [(c \supset a) \supset a]\}$ of the familiar propositional calculus of today. Substituting Peirce's hook or "claw" of illation (—≺) or Schröder's subsumption (∈) for the "horseshoe (⊃) and Peirce's over-bar or Schröder's negation prime for the tilde of negation suffices to yield the same formula in the classical Peirce-Schröder calculus; thus:

> Peano-Russell: $[(\sim c \supset a) \supset (\sim a \supset c)] \supset \{(\sim c \supset a) \supset [(c \supset a) \supset a]\}$
> Peirce: $[(\bar{c} \prec a) \prec (\bar{a} \prec c)] \prec \{(\bar{c} \prec a) \prec [(c \prec a) \prec a]\}$
> Schröder: $[(c' \in a) \in (a' \in c)] \in \{(c' \in a) \in [(c \in a) \in a]\}$

One of Husserl's arguments against the claim that the algebra of logic generally, including the systems of Boole and Peirce, and Schröder's system in particular [Husserl *1891*, 267–272] is the employment of 1 and 0 rather than truth-values *true* and *false*. Certainly neither Boole nor Peirce had not been averse to employing Boolean values (occasionally even using '∞' for universes of discourses of indefinite or infinite cardinality) in analyzing the truth of propositions. Husserl, however, made it a significant condition of his determination of logic as a calculus, as opposed to logic as a language, that truth-values be made manifest, and not represented by numerical values, and he tied this to the mental representation which languages serve.

In the manuscript "On the Algebraic Principles of Formal Logic" written in the autumn of 1879—the very year in which Frege's *Begriffsschrift* appeared, Peirce (see [Peirce *1989*, 23]) explicitly identified his "claw" as the "copula of inclusion" and defined material implication or logical inference, illation, as

> 1st, $A \prec A$, whatever $A$ may be.
> 2nd If $A \prec B$, and $B \prec C$, then $A \prec C$.

From there he immediately connected his definition with truth-functional logic, by asserting [Peirce *1989*, 23] that

> This definition is sufficient for the purposes of formal logic, although it does not distinguish between the relation of inclusion and its converse. Were it desirable thus to distinguish, it would be sufficient to add that the real truth or falsity of $A \prec B$, supposes the existence of $A$.

The following year, Peirce continued along this route: in "The Algebra of Logic" of 1880 [Peirce *1880*, 21; *1989*, 170],

$$A \prec B$$

is explicitly defined as "*A* implies *B*", and

$$A \overline{\prec} B$$

defines "*A* does not imply *B*." Moreover, we are able to distinguish universal and particular propositions, affirmative and negative, according to the following scheme:

> A.  $a \prec b$    All *A* are *B*    (universal affirmative)
> E.  $a \prec \bar{b}$    No *A* is *B*    (universal negative)

| I. | $\breve{a} \prec b$ | Some *A* is *B* | (particular affirmative) |
| O. | $\breve{a} \prec \bar{b}$ | Some *A* is not *B* | (particular negative) |

In 1883 and 1884, in preparing preliminary versions for his article "On the Algebra of Logic: A Contribution to the Philosophy of Notation" [Peirce *1885*], Peirce develops in increasing detail the truth functional analysis of the conditional and presents what we would today recognize as the indirect or abbreviated truth table.

In the undated manuscript "Chapter III. Development of the Notation" [Peirce *n.d.(c)*], composed *circa* 1883-1884, Peirce undertook an explanation of material implication (without, however, explicitly terming it such), and making it evident that what he has in mind is what we would today recognize as propositional logic, asserting that letters represent assertions, and exploring the conditions in which inferences are valid or not, i.e., undertaking to "develope [*sic*] a calculus by which, from certain assertions assumed as premises, we are to deduce others, as conclusions." He explains, further, that we need to know, given the truth of one assertion, how to determine the truth of the other.

And in 1885, in "On the Algebra of Logic: A Contribution to the Philosophy of Notation" [Peirce *1885*], Peirce sought to redefine categoricals as hypotheticals and presented a propositional logic, which he called *icon of the first kind*. Here, Peirce [*1885*, 188–190], Peirce considered the *consequentia*, and introduces inference rules, in particular *modus ponens*, the "icon of the second kind" [Peirce *1885*, 188], transitivity of the copula or "icon of the third kind" [Peirce *1885*, 188–189], and *modus tollens*, or "icon of the fourth kind" [Peirce *1885*, 189].

In the manuscript fragment "Algebra of Logic (Second Paper)" written in the summer of 1884, Peirce (see [Peirce *1986*, 111–115]) reiterated his definition of 1880, and explained in greater detail there [Peirce *1986*, 112] that: "In order to say "If it is *a* it is *b*," let us write $a \prec b$. The formulae relating to the symbol $\prec$ constitute what I have called the algebra of the copula…. The proposition $a \prec b$ is to be understood as true if either *a* is false or *b* is true, and is only false if *a* is true while *b* is false."

It was at this stage that Peirce undertook the truth-functional analysis of propositions and of proofs, and also introduced specific truth-functional considerations, saying that, for **v** is the symbol for "true" (*verum*) and **f** the symbol for false (*falsum*), the propositions $\mathbf{f} \prec a$ and $a \prec \mathbf{v}$ are true, and either one or the other of $\mathbf{v} \prec a$ or $a \prec \mathbf{f}$ are true, depending upon the truth or falsity of *a*, and going on to further analyze the truth-functional properties of the "claw" or "hook".

In Peirce's conception, as found in his "Description of a Notation for the Logic of Relatives" of 1870, then Aristotelian syllogism becomes a hypothetical proposition, with material implication as its main connective; he writes [Peirce *1870*, 518] *Barbara* as

If $x \prec y$,
and $y \prec z$,
then $x \prec z$.

In Frege's *Begriffsschrift* notation of 1879 [Frege *1879*, §6], this same argument would be rendered as:

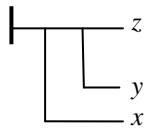

In the familiar Peano-Russell notation, this is just

$$(x \supset y) \cdot (y \supset z)] \supset (x \supset y).$$

Schröder, ironically, even complained about what he took to be Peirce's (and Hugh MacColl's) efforts to base logic on the propositional calculus, which he called the "MacColl-Peircean propositional logic."

Frege [*1895*, 434] recognized that implication was central to the logical systems of Peirce and Schröder (who employed '∈, or Subsumption, in lieu of Peirce's '$\prec$'), although criticizing them for employing the same symbol for class inclusion (or ordering) and implication, and thus for allegedly failing distinguish between these; class and set are in Schröder, he says [Frege *1895*, 435] "eingemischt", and which, in actuality, is just the part-whole relation. Thus he writes [Frege *1895*, 434]:

> Was Herr Schröder ‚Einordnung' oder ‚Subsumption' nennt, ist hier eigentlich nichts Anderes als die Beziehung des Teiles zum Ganzen mit der Erweiterung, dass jedes Ganze als seiner selbst betrachtet werden soll.

Frege [*1895*, 441–442] thus wants Schröder to distinguish between the s*ubter*-Beziehung, the class-inclusion relation, which is effectively implication, referencing [Frege *1895*, 442*n*.] in this regard Peano's [*1894*, §6] '⊃', and the *sub*-Beziehung, or set membership relation, referencing (at [Frege *1895*, 442*n*.]) Peano's [*1894*, §6] '∈'.

John Shosky [*1997*] distinguished between the truth-table *technique* or method on the one hand and the truth-table *device* on the other, the former amounting to a truth-functional analysis of propositions or arguments, the latter being in essence the presentation of truth-functional analysis in a tabular, or matrix, array. On this basis he argued that truth tables first appeared on paper in recognizable form around 1912, composed in the hand of either Ludwig Wittgenstein, with an addition by Bertrand Russell, on the verso of a typescript of a lecture by Russell on logical atomism, and thus earlier than its appearance in Wittgenstein's *Tractatus Logico-philosophicus* [Wittgenstein *1922*, Prop. 4.31] or the work of Emil Leon Post (1897–1954) and Jan Łukasiewicz (1878–1956) in 1920.[1]

The first instance by Peirce of a truth-functional analysis which satisfies the conditions for truth tables, but is not yet constructed in tabular form, is in his 1885 article "On the Algebra of Logic: A Contribution to the Philosophy of Notation", in which he gave a proof, using the truth table method, of what has come to be known as *Peirce Law*: $((A \rightarrow B) \rightarrow A) \rightarrow A$, his "fifth icon", whose validity he tested using truth-functional analysis. In an untitled paper written in 1902 and placed in volume 4 of the Hartshorne and Weiss edition of Peirce's *Collected Papers*, Peirce displayed the following table for three terms, *x*, *y*, *z*, writing **v** for *true* and **f** for *false* ("The Simplest Mathematics"; January 1902 ("Chapter III. The Simplest Mathematics (Logic III))", RC MS #431, January 1902; see [Peirce *1933*, 4:260–262]).

| *x* | *y* | *z* |
|---|---|---|
| **v** | **v** | **v** |
| **v** | **f** | **f** |
| **f** | **v** | **f** |
| **f** | **f** | **v** |

where *z* is the term derived from an [undefined] logical operation on the terms *x* and *y*.[2] In February 1909, while working on his trivalent logic, Peirce applied the tablular method to various connectives, for example negation of *x*, as $\bar{x}$, written out, in his notebook ([Peirce *1865-1909*]; see [Fisch & Turquette *1966*]), as:[3]

---

[1] [Shosky *1997*, 12, *n*. 6] cites Post's [*1921*] in [van Heijenoort *1967a*, 264–283]; but see also the doctoral thesis [Post *1920*], which he misses. Łukasiewicz is mentioned, but [Shosky *1997*] gives no reference; see [Łukasiewicz *1920*]. See also [Anellis *2004*] for historical details.

[2] [Shosky *1997*] totally ignores the detailed and complex history of the truth-table method and shows no knowledge of the existence of the truth-table device of Peirce in 1902. For a critique of Shosky's [*1997*] account and the historical background to Peirce's work, see [Anellis *2004*]; see also [Clark *1997*] and [Zellweger *1997*] for their re-(dis)-covery and exposition of Peirce's work on truth-functional analysis and the development of his truth-functional matrices.

[Grattan-Guinness *2004-05*, 187–188], meanwhile, totally misrepresents the account in [Anellis *2004*] of the contributions of Peirce and his cohorts to the evolution and development of truth-functional analysis and truth tables in suggesting that: (1) Peirce did not achieve truth table matrices and (2) that [Anellis *2004a*] was in some manner attempting to suggest that Russell somehow got the idea of truth tables from Peirce. The latter is actually strongly contraindicated on the basis of the evidence that was provided in [Anellis *2004*], where it is shown that Peirce's matrices appeared in unpublished manuscripts which did not arrive at Harvard until the start of 1915, after Russell had departed, and were likely at that point unknown, so that, *even if* Russell could had been made aware of them, it would have more than likely have been from Maurice Henry Sheffer (1882–1964), and *after* 1914.

[3] Peirce's orginal tables from MS 399 are reproduced as plates 1-3 at [Fisch & Turquette *1966*, 73–75].

| $x$ | $\bar{x}$ |
|---|---|
| **V** | **F** |
| **L** | **L** |
| **F** | **V** |

where, **V**, **F**, and **L** are the truth-values true, false, and indeterminate or unknown respectively, which he called "limit".[4] Russell's rendition of Wittgenstein's tabular definition of negation, as written out on the verso of a page from Russell's transcript notes, using 'W' (wahr) and 'F' (falsch) (see [Shosky 1997, 20]), where the negation of *p* is written out by Wittgenstein as p ◊ q, with Russell adding "= ~p", to yield: p ◊ q = ~p is

| p | q |
|---|---|
| W | W |
| W | F |
| F | W |
| F | F |

The trivalent equivalents of classical disjunction and conjunction were rendered by Peirce in that manuscript respectively as

| Θ | V | L | F |     | Z | V | L | F |
|---|---|---|---|-----|---|---|---|---|
| **V** | V | V | V |     | **V** | V | L | F |
| **L** | V | L | L |     | **L** | L | L | F |
| **F** | V | L | F |     | **F** | F | F | F |

Max Fisch and Atwell Turquette [*1966*, 72], referring to [Turquette *1964*, 95–96], assert that the tables for trivalent logic in fact were extensions of Peirce's truth tables for bivalent logic, and hence prior to 23 February 1909 when he undertook to apply matrices for the truth-functional analysis for trivalent logic. The reference is to that part of Peirce's [*1885*, 183–193], "On the Algebra of Logic: A Contribution to the Philosophy of Notation"—§II "Non-relative Logic"—dealing with truth-functional analysis, and Turquette [*1964*, 95] uses "truth-function analysis" and "truth-table" synonymously, a confusion which, in another context, [Shosky *1997*] when warning against confusing, and insisting upon a careful distinction between the truth-table *technique* and the truth-table *device*.

    Roughly contemporary with the manuscript "The Simplest Mathematics" is "Logical Tracts. No. 2. On Existential Graphs, Euler's Diagrams, and Logical Algebra", *ca*. 1903 [Peirce *1933*, 4.476]; Harvard Lectures on Pragmatism, 1903 [Peirce *1934*, 5.108]),

    In the undated manuscript [Peirce *n.d.(b)*] identified as composed *circa* 1883-84 "On the Algebra of Logic" and the accompanying supplement, we find what unequivocally would today be labeled as an indirect or abbreviated truth table for the formula {((a —≺b)—≺c)—≺d}—≺e, as follows:

$$\overline{\{((a \prec b) \prec c) \prec d\} \prec e}$$
$$\begin{array}{cccccc} f & f & f & f & \prec f \\ f & v & & v & f \\ - & - & - & - & v \end{array}$$

    The whole of the undated eighteen-page manuscript "Logic of Relatives", also identified as composed *circa* 1883-84 [Peirce *n.d.(c)*; MS #547], is devoted to a truth-functional analysis of the conditional, which includes the equivalent, in list form, of the truth table for $x \prec y$, as follows [Peirce *n.d.(c)*; MS #547:00016; 00017]:

---

    [4] In *Logic (Logic Notebook 1865–1909)*; MS 339:440–441, 443; see [Peirce *1865–1909*], which had been examined by Fisch and Turquette. On the basis of this work, Fisch and Turquette [*1966*, p. 72] concluded that by 23 February 1909 Peirce was able to extend his truth-theoretic matrices to three-valued logic, thereby anticipating both Jan Łukasiewicz in "O logice trójwartosciowej" [Łukasiewicz *1921*], and Emil Leon Post in "Introduction to a General Theory of Elementary Propositions" (Post *1921*), by a decade in developing the truth table device for triadic logic and multiple-valued logics respectively. Peirce's tables from MS 399 are reproduced as the plates at [Fisch & Turquette *1966*, 73–75].

$$x \prec y$$

|is true when|is false when|
|---|---|
|$x = \mathbf{f}\ \ y = \mathbf{f}$|$x = \mathbf{v}\ \ y = \mathbf{f}$|
|$x = \mathbf{f}\ \ y = \mathbf{v}$| |
|$x = \mathbf{v}\ \ y = \mathbf{v}$| |

Peirce also wrote follows [Peirce *n.d.(c)*; MS #547:00016] that: "It is plain that $x \prec y \prec z$ is false only if x = **v**, (*y* $\prec z$) = **f**, that is only if $x = \mathbf{v}$, $y = \mathbf{v}$, $z = \mathbf{f}$…."

Finally, in the undated manuscript "An Outline Sketch of Synechistic Philosophy" identified as composed in 1893, we have an unmistakable example of a truth table matrix for a proposition and its negation [Peirce *1893*; MS #946:00004], as follows:

|   | t | f |
|---|---|---|
| t | t | f |
| f | t | t |

which is clearly and unmistakably equivalent to the truth-table matrix for $x \prec y$ in the contemporary configuration, expressing the same values as we note in Peirce's list in the 1883-84 manuscript "Logic of Relatives" [Peirce *n.d.(c)*; MS #547:00016; 00017]. That the multiplication matrices are the most probable inspiration for Peirce's truth-table matrix is that it appears alongside matrices for a multiplicative two-term expression of linear algebra for $\{i, j\}$ and $\{i, i - j\}$ [Peirce *1893*; MS #946:00004]. Indeed, it is virtually the same table, and in roughly—*i.e.*, apart from inverting the location within the respective tables for antecedent and consequent—the same configuration as that found in the notes, taken in April 1914 by Thomas Stearns Eliot (1888–1965) in Russell's Harvard University logic course (as reproduced at [Shosky *1997*, 23]), where we have:

$p \vee q$    $q\begin{cases}\begin{array}{c|cc} & T & F \\ \hline T & T & T \\ F & T & F \end{array}\end{cases}$    $p \supset q$    $q\begin{cases}\begin{array}{c|cc} & T & F \\ \hline T & T & T \\ F & F & T \end{array}\end{cases}$    $\sim p \vee \sim q$    $q\begin{cases}\begin{array}{c|cc} & T & F \\ \hline T & T & T \\ F & T & F \end{array}\end{cases}$

Finally, and also contemporaneous with this work, and continuing to experiment with notations, Peirce developed his "box-X" or "X-frame" notation, which resemble the square of opposition in exploring the relation between the relations between two terms or propositions. Lines may define the perimeter of the square as well as the diagonals between the vertices; the presence of connecting lines between the corners of a square or box indicates the states in which those relations are false, or invalid, absence of such a connecting line indicates relations in which true or valid relation holds. In particular, as part of this work, Peirce developed a special iconic notation for the sixteen binary connectives, as follows (from "The Simplest Mathematics" written in January 1902 ("Chapter III. The Simplest Mathematics (Logic III)", MS 431; [Peirce *1902*]; see [Clark *1997*, 309]) containing a table presenting the 16 possible sets of truth values for a two-term proposition:

| 1 | 2 | 3 | 4 | 5 | 6 | 7 | 8 | 9 | 10 | 11 | 12 | 13 | 14 | 15 | 16 |
|---|---|---|---|---|---|---|---|---|----|----|----|----|----|----|----|
| F | F | F | F | T | T | T | T | F | F  | F  | F  | T  | T  | T  | T  |
| F | F | F | T | F | T | F | F | T | T  | F  | T  | F  | T  | T  | T  |
| F | F | T | F | F | F | T | F | F | T  | T  | T  | T  | F  | T  | T  |
| F | T | F | F | F | F | F | T | F | T  | T  | T  | T  | T  | F  | T  |

that enabled him to give a quasi-mechanical procedure for identifying thousands of tautologies from the substitution sets of expressions with up to three term (or proposition) variables and five connective variables.

In his X-frames notation, the open and closed quadrants are indicate truth or falsity respectively, so that for example, ⊠, the completely closed frame, represents row 1 of the table for the sixteen binary connectives, in which all assignments are false, and x, the completely open frame, represents row 16, in which all values are true (for details, see [Clark *1997*] and [Zellweger *1997*]). The X-frame notation is based on the representation of truth-values for two terms as follows:

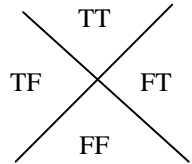

The full details of this scheme are elaborated by Peirce in his manuscript "A Proposed Logical Notation (Notation)" of *circa* 1903 [Peirce *ca. 1903*, esp. 530:0000026-28].

  **3. Conclusion.** In a manuscript of 1893, in the context of his study of the truth-functional analysis of propositions and proofs and his continuing efforts at defining and understanding the nature of logical inference, and against the background of his mathematical work in matrix theory in algebra, Charles Peirce presented a truth table which displayed in matrix form the definition of his most fundamental connective, that of illation, which is equivalent to the truth-functional definition of material implication. Peirce's matrix is exactly equivalent to that for material implication discovered by Shosky that is attributable to Bertrand Russell and has been dated as originating in 1912. Thus, Peirce's table of 1893 may be considered to be the earliest known instance of a truth table device in the familiar form which is attributable to an identifiable author, and antedates not only the tables of Post, Wittgenstein, and Łukasiewicz of 1920-22, but Russell's table of 1912 and also Peirce's previously identified tables for trivalent logic tracable to 1902.

*Peirce Edition Project, Institute for American Thought*
*Indiana University – Purdue University at Indianapolis*
ianellis@iupui.edu